\documentclass[MSNbibl,number,citesort,dvips]{arxbj}
\usepackage{upgreek}
\usepackage{graphicx}

\aid{0}
\volume{20}
\issue{2}
\pubyear{2014}
\firstpage{1006}
\lastpage{1028}
\doi{10.3150/13-BEJ514} 

\makeatletter
\renewcommand{\pi}{\uppi}
\newcommand{\llvert}{\vert}
\def\nosymbol{}
\def\sign{\operatorname{sign}}%
\newtheorem{theorem}{Theorem}
\newtheorem{lemma}{Lemma}

\newcommand{\cov}{\operatorname{cov}}

\def\sign{\operatorname{sign}}
\def\cvm{Cram{\'e}r von Mises}
\def\dcov{\operatorname{dcov}}
\def\mR{\mathbb{R}}

\newcommand{\tmop}[1]{\mathrm{#1}}

\makeatother

\begin{document}
\begin{frontmatter}

\title{A consistent test of independence based on a~sign covariance
related to Kendall's tau}
\runtitle{A consistent test of independence}

\begin{aug}
\author{\fnms{Wicher} \snm{Bergsma}\corref{}\thanksref{e1}\ead[label=e1,mark]{w.p.bergsma@lse.ac.uk}}
\and
\author{\fnms{Angelos} \snm{Dassios}\thanksref{e2}\ead[label=e2,mark]{a.dassios@lse.ac.uk}}
\runauthor{W. Bergsma and A. Dassios} 
\address{London School of Economics and Political Science, Houghton
Street, London WC2A 2AE, United
Kingdom. \printead{e1,e2}}
\end{aug}

\received{\smonth{7} \syear{2012}}
\revised{\smonth{1} \syear{2013}}

%
\begin{abstract}
The most popular ways to test for independence of two ordinal random
variables are by means of Kendall's tau and Spearman's rho. However,
such tests are not consistent, only having power for alternatives with
``monotonic'' association. In this paper, we introduce a natural
extension of Kendall's tau, called $\tau^*$, which is non-negative and
zero if and only if independence holds, thus leading to a consistent
independence test. Furthermore, normalization gives a rank correlation
which can be used as a measure of dependence, taking values between
zero and one.
A comparison with alternative measures of dependence for ordinal random
variables is given, and it is shown that, in a well-defined sense,
$\tau
^*$ is the simplest, similarly to Kendall's tau being the simplest of
ordinal measures of monotone association. Simulation studies show our
test compares well with the alternatives in terms of average $p$-values.
\end{abstract}

%
\begin{keyword}
\kwd{concordance}
\kwd{copula}
\kwd{discordance}
\kwd{measure of association}
\kwd{ordinal data}
\kwd{permutation test}
\kwd{sign test}
\end{keyword}

\end{frontmatter}

\section{Introduction}

A random variable $X$ is called \emph{ordinal} if its possible values
have an ordering, but no
distance is assigned to pairs of outcomes. Ordinal variables may be
continuous, categorical, or
mixed continuous/categorical.
Ordinal data frequently arise in many fields, though especially often
in social and biomedical
science (Kendall and Gibbons \cite{kg90}, Agresti \cite{agresti10}). Ordinal data methods are also often
applied to real-valued
(interval level) data in order to achieve robustness.

\begin{figure}
\centering
\begin{tabular}{@{}cc@{}}

\includegraphics{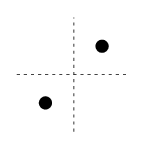}
 & \includegraphics{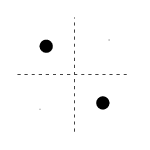}\\
\footnotesize{(a) Concordant pair} & \footnotesize{(b) Discordant pair }
\end{tabular}
\caption{Concordant and discordant pairs of points
associated with Kendall's tau.}\label{concdisc2}
\end{figure}

\begin{figure}[b]
\centering
\begin{tabular}{@{}cccc@{}}

\includegraphics{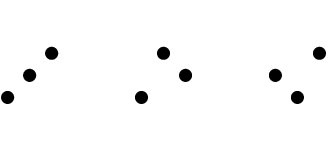}
 &&& \includegraphics{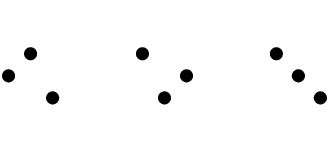}\\
\footnotesize{(a) Concordant triples} &&& \footnotesize{(b) Discordant triples}\\
\end{tabular}
\caption{Concordant and discordant triples of points
associated with Spearman's rho.}\label{concdisc3}
\end{figure}

The two most popular measures of association for ordinal random
variables $X$ and $Y$ are Kendall's tau ($\tau$) (Kendall \cite{kendall38}) and
Spearman's rho ($\rho_S$) (Spearman \cite{spearman04}), which may be defined as
\[
\tau= E\sign\bigl[(X_1-X_2) (Y_1-Y_2)
\bigr],\qquad  \rho_S=3E\sign \bigl[(X_1-X_2)
(Y_1-Y_3)\bigr],
\]
where the $(X_i,Y_i)$ are independent replications of $(X,Y)$ (Kruskal \cite{kruskal58}).
The factor 3 in the expression for $\rho_S$ occurs to obtain a measure
whose range is $[-1,1]$.
Both $\tau$ and $\rho_S$ are proportional to sign versions of the
ordinary covariance, which can be seen from the following expressions
for the covariance:
\[
\cov(X,Y) = \tfrac{1}{2}E(X_1-X_2)
(Y_1-Y_2) = E(X_1-X_2)
(Y_1-Y_3).
\]
From the definitions, probabilistic interpretations of $\tau$ and
$\rho
_S$ can be derived. Firstly,
%
\begin{equation}
\label{ktauprob} \tau= \Pi_{C_2} - \Pi_{D_2},
\end{equation}
where $\Pi_{C_2}$ is the probability that two observations are
concordant and $\Pi_{D_2}$ the probability that they are discordant
(see Figure~\ref{concdisc2}). Secondly,
\[
\rho_S = \Pi_{C_3} - \Pi_{D_3},
\]
where $\Pi_{C_3}$ is the probability that three observations are
concordant and $\Pi_{D_3}$ the probability that they are discordant
(see Figure~\ref{concdisc3}).
It can be seen that $\tau$ is simpler than $\rho_S$, in the sense that
it can be defined using only two rather than three independent
replications of $(X,Y)$, or, more specifically, in terms of
probabilities of concordance and discordance of two rather than three
points. This was a reason for Kruskal to prefer $\tau$ to $\rho_S$
(Kruskal \cite{kruskal58}, end of Section 14).

An alternative definition of $\rho_S$, which was originally given by
Spearman, is as a Pearson correlation between uniform rank scores of
the $X$ and $Y$ variables. For continuous random variables, both this
and the aforementioned definition lead to the same quantity. However,
with this definition, $\rho_S$ is to some extent an \emph{ad hoc}
measure, since the choice of scores is arbitrary, and alternative
scores (e.g., normal scores) might be used.

A test of independence based on i.i.d. data can be obtained by
application of the permutation test to an estimator of $\tau$ or $\rho
_S$, which is easy to implement and fast to carry out with modern computers.
Such ordinal tests are also used as a robust alternative to tests based
on the Pearson correlation.

A drawback for certain applications is that $\tau$ and $\rho_S$ may be
zero even if there is
an association between $X$ and $Y$, so tests based on them are inconsistent
{{for the alternative of a general association}}. For this reason, alternative
coefficients have been devised. The best known of these are those
introduced by
Hoeffding \cite{hoeffding48ind} and Blum, Kiefer and Rosenblatt \cite{bkr61}.
With $F_{12}$ the joint distribution function of $(X,Y)$, and $F_1$ and
$F_2$ the marginal distribution functions of $X$, respectively, $Y$,
Hoeffding's coefficient is given as
%
\begin{equation}
\label{hoeff} H = \int\bigl[F_{12}(x,y)-F_1(x)F_2(y)
\bigr]^2\,\mathrm{d}F_{12}(x,y),
\end{equation}
and the Blum--Kiefer--Rosenblatt (henceforth: BKR) coefficient as
%
\begin{equation}
\label{dd} D = \int\bigl[F_{12}(x,y)-F_1(x)F_2(y)
\bigr]^2\,\mathrm{d}F_{1}(x)\,\mathrm{d}F_2(y).
\end{equation}
Both can be seen to be non-negative with equality to zero under
independence. Furthermore, $D=0$ can also be shown to imply independence.
However, the Hoeffding coefficient has a severe drawback, namely that
it may be zero even if there is an association, that is, it does not
lead to a consistent independence test. An example is the case that
$P(X=0,Y=1)=P(X=1,Y=0)=1/2$ (Hoeffding \cite{hoeffding48ind}, page 548).

A third option, especially suitable for categorical data, is the
Pearson chi-square test; it is directly applicable to categorical data
and can be used for continuous data after a suitable categorization.
However, the chi-square test does not take the ordinal nature of the
data into account, leading to potential power loss for ``ordinal''
alternatives; effectively the chi-square test treats the data as
nominal rather than ordinal (see also Agresti \cite{agresti10}).

Although $H$ and $D$ have simple mathematical formulas, they seem to be
rather arbitrary, and many variants are possible (see also Section~\ref
{sec-rel}).
For this reason, we decided to develop a probabilistic interpretation
of $H$ (given in Section~\ref{sec-alt} of this paper). However, we then
noticed that $H$ and $D$ were unnecessarily complex, and that a clearly
simpler and natural alternative coefficient was possible. Our new
coefficient is a direct modification of Kendall's $\tau$, which we call
$\tau^*$. It is non-negative and zero if and only if independence holds.
Like $\tau$ and $\rho_S$, we show that $H$ and $\tau^*$ equal the
difference of concordance and discordance probabilities of a number of
independent replications of $(X,Y)$.
Analogously to the aforementioned way that $\tau$ is simpler than
$\rho
_S$, $\tau^*$ is simpler than $H$ in that only four independent
replications of $(X,Y)$ are required, whereas $H$ needs five.
It appears to us that relative simplicity of interpretation of a
coefficient is of utmost importance, and that this is also the main
reason for the current popularity of Kendall's tau. In particular, when
it was introduced in the pre-computer age in 1938, the sample value of
Kendall's tau was much harder to compute than the sample value of
Spearman's rho, which had been in use since 1904 (Kruskal \cite{kruskal58}). In
spite of this, judging by the number of Google Scholar hits, both
currently appear to be about equally popular.\footnote{The Google
Scholar search ``\texttt{kendall's tau}'' \texttt{OR} ``\texttt{kendall tau}'' gave us 16,400
hits and the search ``s\texttt{pearman's rho}'' \texttt{OR} ``\texttt{spearman rho}'' 18,500.}

As a remark on the two-sample case, if one of the variables is binary,
a test that $\tau^*=0$ is equivalent to the \cvm\ test, as shown in
Section~3 in Dassios and Bergsma~\cite{db12}.

The organization of this paper is as follows.
In Section~\ref{sec2}, we first define $\tau^*$, and then state our
main theorem that $\tau^*\ge0$ with equality if and only if
independence holds. Furthermore, we provide a probabilistic
interpretation in terms of concordance and discordance probabilities of
four points.
Section~\ref{sec3} contains the proof of the main theorem. The proof
turns out to be surprisingly involved for such a simple to formulate
coefficient, and the ideas in the proof may be useful for other related
research.
A comparison with the Hoeffding, the BKR and some more recent
coefficients is given in Section~\ref{sec-alt}, and a new probabilistic
interpretation for the former is given.
In Section~\ref{sec-applic}, we give a description of independence
testing via the permutation test and a simulation study compares
average $p$-values of our test and the aforementioned other two tests.
Our test compares well with the other two in this respect.


\section{\texorpdfstring{Definition of $\tau^*$ and statement of its properties}
{Definition of tau* and statement of its properties}}\label{sec2}

We denote i.i.d. sample values by $(x_1,y_1),\ldots,(x_n,y_n)$, but
will also use $\{(X_i,Y_i)\}$ to denote i.i.d. replications of $(X,Y)$
in order to define population coefficients.
The empirical value $t$ of Kendall's tau is
\[
t = \frac{1}{n^2}{\sum_{i,j=1}^n
\sign(x_i-x_j)\sign(y_i-y_j)},
\]
and its population version is
\[
\tau= {E\sign(X_1-X_2)\sign(Y_1-Y_2)}.
\]
%
(Kruskal \cite{kruskal58}, Kendall and Gibbons \cite{kg90}). With
\begin{eqnarray*}
s(z_1,z_2,z_3,z_4) &=& {
\sign(z_1-z_4) (z_3-z_2)}
\\
&=& \sign\bigl(|z_1-z_2|^2+|z_3-z_4|^2-|z_1-z_3|^2-|z_2-z_4|^2
\bigr),
\end{eqnarray*}
we obtain
\[
t^2 = \frac{1}{n^4}{\sum_{i,j,k,l=1}^n
s(x_i,x_j,x_k,x_l)s(y_i,y_j,y_k,y_l)}
\]
and
\begin{eqnarray*}
\tau^2 = {Es(X_1,X_2,X_3,X_4)s(Y_1,Y_2,Y_3,Y_4)}.
\label{tdef}
\end{eqnarray*}
Replacing squared differences in $s$ by absolute values of differences,
we define
%
\begin{equation}
a(z_1,z_2,z_3,z_4) = \sign
\bigl(|z_1-z_2|+|z_3-z_4|-|z_1-z_3|-|z_2-z_4|
\bigr). \label{sdef}
\end{equation}
This leads to a modified version of $t^2$,
%
\begin{equation}
t^* = \frac{1}{n^4}\sum_{i,j,k,l=1}^na(x_i,x_j,x_k,x_l)a(y_i,y_j,y_k,y_l)
\label{tdef2}
\end{equation}
and the corresponding population coefficient
\[
\tau^* = \tau^*(X,Y) = {Ea(X_1,X_2,X_3,X_4)a(Y_1,Y_2,Y_3,Y_4)}.
\label{tdef3}
\]
The quantities $t^*$ and $\tau^*$ are new, and the main result of the
paper is the following:
%
\begin{theorem}\label{th1}
Assume $(X,Y)$ has a bivariate discrete or continuous distribution, or
a mixture of the two, that is, assume there exists a probability mass
function $f$ and a density function $\tilde{f}$ such that
\[
P ( X<x, Y<y )= \sum_{u_i < x, v_i <y } f ( u_i,
v_i ) + \int_{u < x, v < y} \tilde{f} ( u,v ) \,\tmop{d}u
\,\tmop{d}v.
\]
It holds true that $\tau^*(X,Y)\ge0$ with equality if and only if $X$
and $Y$ are independent.
\end{theorem}
%
The proof is given in Section~\ref{sec3}.
We conjecture that the condition of the theorem is not necessary, that
is, that for arbitrary $(X,Y)$ it holds that $\tau^*(X,Y)\ge0$ with
equality if and only if $X$ and~$Y$ are independent.

If the sign functions are omitted from $\tau^*$, we obtain the
covariance introduced by Bergsma \cite{bergsma06} and Sz{\'e}kely, Rizzo and Bakirov \cite{srb07}. They
showed that for arbitrary real random variables~$X$ and~$Y$, this
covariance is non-negative with equality to zero if and only if~$X$ and~$Y$ are independent. (See Section~\ref{non-ord} for further details.)


By the Cauchy--Schwarz inequality, the normalized value
\[
\tau_b^* = \frac{\tau^*(X,Y)}{\sqrt{\tau^*(X,X)\tau^*(Y,Y)}}
\]
does not exceed one. (Note that this notation is in line with Kendall's
$\tau_b$, defined analogously.)

The definition of $\tau^*$ can easily be extended to $X$ and $Y$ in
arbitrary metric spaces, but unfortunately Theorem~\ref{th1} does not
extend then, as it is possible that $\tau^*<0$. This is shown by the
following example. Consider a set of points $\{u_1,\ldots,u_8\}\subset
\mR^8$, where $u_i=(u_{i1},\ldots,u_{i8})'$ such that $u_{ii}=3$,
$u_{ij}=-1$ if $i\ne j$ and $i,j\le4$ or $i,j\ge5$,\vadjust{\goodbreak} and $u_{ij}=0$
otherwise. Suppose $Y$ is uniformly distributed on $\{0,1\}$, and given
$Y=0$, $X$ is uniformly distributed on $u_1,\ldots,u_4$, and given
$Y=1$, $X$ is uniformly distributed on $u_5,\ldots,u_8$. Then $\tau^*=-1/64$.

Note that $\tau^*(X,Y)$ is a function of the copula, which is the joint
distribution of $F_1(X)$ and $F_2(Y)$, where $F_1$ and $F_2$ are the
cumulative distribution functions of $X$ and $Y$. More generally, for
any strictly monotone (increasing or decreasing) functions $g$ and $h$,
$\tau^*(X,Y)=\tau^*(g(X),h(Y))$. Nelsen \cite{nelsen06}, Chapter~5, explores
the way in which copulas can be used in the study of dependence between
random variables, paying particular attention to Kendall's tau and
Spearman's rho.


We now give a probabilistic interpretation of $\tau^*$.
Recall that Kendall's tau is the probability that a pair of points is
concordant minus the probability that a pair of points is discordant.
Our $\tau^*$ is proportional to the probability that two pairs are
``jointly'' concordant, plus the probability that two pairs are ``jointly''
discordant, minus the probability that, ``jointly'', one pair is
discordant and the other concordant. Here, ``jointly'' refers to there
being a common axis separating the two points of each of the two pairs.

\begin{figure}
\centering
\begin{tabular}{@{}ccccc@{}}

\includegraphics{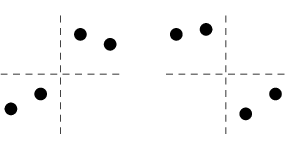}
 &&&& \includegraphics{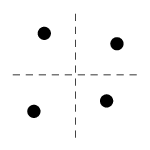}\\
\footnotesize{(a) Concordant points} &&&& \footnotesize{(b) Discordant points}
\end{tabular}
\caption{Configurations of concordant and discordant quadruples of
points associated
with $\tau^*$. The dotted axes indicate strict separation of points in
different quadrants; within a quadrant, no restrictions apply on the
relative positions of points.}\label{concdisc}
\end{figure}

To use a slightly different terminology which will be convenient, we
say that a set of four points is concordant if two pairs are either
``jointly'' concordant or ``jointly'' discordant, while four points are
called discordant if, ``jointly'', one pair is concordant and the other
is discordant. These configurations are given in Figure~\ref{concdisc}.
In mathematical notation, a~set of four points $\{(x_1,y_1),\ldots
,(x_4,y_4)\}$
is concordant if there is a permutation $(i,j,k,l)$ of $(1,2,3,4)$ such that
\[
(x_i,x_j<x_k,x_l) \&
\bigl[(y_i,y_j<y_k,y_l)||
(y_i,y_j>y_k,y_l) \bigr],
\]
and discordant if there is a permutation $(i,j,k,l)$ of $(1,2,3,4)$
such that
\begin{eqnarray*}
\bigl[(x_i,x_j<x_k,x_l)||(x_i,x_j>x_k,x_l)
\bigr] \& \bigl[(y_i,y_k<y_j,y_l)||
(y_i,y_k>y_j,y_l) \bigr],
\end{eqnarray*}
where $||$ and $\&$ are logical OR, respectively, AND, and $I(z_1,z_2<z_3,z_4)$
is shorthand for $I(z_1<z_3\& z_1<z_4\& z_2<z_3\& z_2<z_4)$.
It is straightforward to verify that
\begin{eqnarray*}
a(z_1,z_2,z_3,z_4) &=&
{I(z_1,z_3<z_2,z_4) +
I(z_1,z_3>z_2,z_4)}
\\
&& {}- I(z_1,z_2<z_3,z_4) -
I(z_3,z_4<z_1,z_2),
\end{eqnarray*}
where $I$ is the indicator function. Hence,
%
\begin{eqnarray}\label{tauprob1}
\tau^* &=& {4}P(X_1,X_2<X_3,X_4
\& Y_1,Y_2<Y_3,Y_4)
\nonumber
\\
&&{}+ {4}P(X_1,X_2<X_3,X_4\&
Y_1,Y_2>Y_3,Y_4)
\\
&&{}  -{8}P(X_1,X_2<X_3,X_4\&
Y_1,Y_3<Y_2,Y_4).
\nonumber
\end{eqnarray}
Denoting the probability that four randomly chosen points are
concordant as $\Pi_{C_4}$ and the probability that they are discordant as
$\Pi_{D_4}$, we obtain that the sum of the first two terms on the right-hand side of~(\ref{tauprob1}) equals $\Pi_{C_4}/6$, while the last term
equals $\Pi_{D_4}/24$.
Hence,
%
\begin{equation}
\tau^* = \frac{2\Pi_{C_4} - \Pi_{D_4}}{{3}}. \label{tau-int}
\end{equation}
%
It can be seen that $t^*$ and $\tau^*$ do not depend on the scale at
which the variables are measured, but only on the ranks or grades of
the observations.
Four points are said to be \emph{tied} if they are neither concordant
nor discordant.
Clearly, for continuous distributions the probability of tied
observations is zero. Hence, under independence, when all
configurations are
equally likely, $\Pi_{C_4}=1/3$ and $\Pi_{D_4}=2/3$, and if one
variable is a strictly monotone function of the other, then
$\Pi_{C_4}=1$ and $\Pi_{D_4}=0$.

\section{Comparison to other tests}\label{sec-alt}

The two most popular (almost) consistent tests of independence for
ordinal random variables are those based on Hoeffding's $H$ and BKR's
$D$, given in~(\ref{hoeff}) and~(\ref{dd}).
We compare $\tau^*$ with these coefficients as well as with the
recently introduced non-ordinal measures of Sz{\'e}kely, Rizzo and Bakirov \cite{srb07} and Gretton \textit{et~al.} \cite{gbss05}.
We give a probabilistic interpretation for $H$ and show that $\tau^*$
is simpler. Since $H=0$ does not imply independence if the
distributions are discrete, it should perhaps not be used, and we are
left with two ordinal coefficients, $\tau^*$ and $D$, of which $\tau^*$
is the simplest.
Further discussions of ordinal data and non-parametric methods for
independence testing are given Agresti \cite{agresti10}, Hollander and Wolfe \cite{hw99} and Sheskin \cite{sheskin07}.

\subsection{Probabilistic interpretation of Hoeffding's $H$}\label{sec-hoeff}

Hoeffding's \cite{hoeffding48ind} coefficient for measuring
deviation from independence for a bivariate distribution function is
given by~(\ref{hoeff}) (see also Blum, Kiefer and Rosenblatt \cite{bkr61}, Hollander and Wolfe \cite{hw99} and Wilding and Mudholkar \cite{wm08}).
An alternative formulation given by Hoeffding is
\[
H = \tfrac14 E\phi(X_1,X_2,X_3)
\phi(X_1,X_4,X_5)\phi (Y_1,Y_2,Y_3)
\phi (Y_1,Y_4,Y_5),
\]
where $\phi(z_1,z_2,z_3)=I(z_1\ge z_2)-I(z_1\ge z_3)$.
Hoeffding's $H$ can be zero for some discrete dependent $(X,Y)$. An
example is the case that $P(X=0,Y=1)=P(X=1,Y=0)=1/2$ (Hoeffding \cite{hoeffding48ind}, page
548).

Interestingly, Hoeffding's $H$ has an interpretation in terms of
concordance and discordance probabilities closely related to the
interpretation of $\tau^*$.
With
\begin{eqnarray*}
F_{12}(x,y) &=& P(X\le x,Y\le y),
\\
F_{1\overline2}(x,y) &=& P(X\le x,Y> y) = F_1(x)-F_{12}(x,y),
\\
F_{\overline12}(x,y) &=& P(X> x,Y\le y) = F_2(y)-F_{12}(x,y),
\\
F_{\overline1\overline2}(x,y) &=& P(X>x,Y>y) = 1-F_1(x)-F_2(y)+F_{12}(x,y),
\end{eqnarray*}
we have the equality
%
\begin{equation}
F_{12}-F_1F_2 = F_{12}F_{\overline1\overline2}
- F_{1\overline
2}F_{\overline12}. \label{f12}
\end{equation}
Let five points be $H$-concordant if four are configured as in
Figure~\ref{concdisc}(a) and the fifth is on the
point where the axes cross and, analogously, five points are
$H$-discordant if four are configured as in Figure~\ref{concdisc}(b)
and the fifth is on the
point where the axes cross. Denote the probabilities of $H$-concordance
and discordance by $\Pi_{C_5}$ and $\Pi_{D_5}$.
Then, omitting the arguments $x$ and $y$,
\[
\int \bigl(F_{12}^2F_{\overline1\overline2}^2+F_{1\overline
2}^2F_{\overline12}^2
\bigr)\,\mathrm{d}F_{12} = \frac{2!2!1!}{5!}\Pi_{C_5} = \frac1{30}
\Pi_{C_5}
\]
and
\[
\int F_{12}F_{1\overline2}F_{\overline12}F_{\overline1\overline
2}\,\mathrm{d}F_{12}
= \frac1{5!}\Pi_{D_5} = \frac1{120}\Pi_{D_5}.
\]
Hence, using~(\ref{f12}),
\[
H = \int (F_{12}F_{\overline1\overline2} - F_{1\overline
2}F_{\overline12}
)^2\,\mathrm{d}F_{12} = \frac{2\Pi_{C_5}-\Pi_{D_5}}{60}.
\]

We can see that Hoeffding's $H$ has two drawbacks compared to $\tau^*$.
Firstly, it is more complex in that it is based on concordance and
discordance of five points rather than four and, secondly, it can be
zero under dependence for certain discrete distributions.

\subsection{The Blum--Kiefer--Rosenblatt coefficient and Spearman's rho}

The coefficient $D$ is given by~(\ref{dd}), and tests based on it were
first studied by Blum, Kiefer and Rosenblatt \cite{bkr61}.
It follows from results in Bergsma \cite{bergsma06} that in the continuous
case, with
%
\begin{eqnarray}
\label{hdef} h(z_1,z_2,z_3,z_4)
&=& |z_1-z_2|+|z_3-z_4|-|z_1-z_3|-|z_2-z_4|,
\nonumber
\\
D &=& Eh\bigl(F_1(X_1),F_1(X_2),F_1(X_3),F_1(X_4)
\bigr)\\
&&{}\times h\bigl(F_2(Y_1),F_2(Y_2),F_2(Y_3),F_2(Y_4)
\bigr).\nonumber
\end{eqnarray}
%
A similar formulation was given by Feuerverger \cite{feuerverger93}, who used
characteristic functions for its derivation. This connection of
Feuerverger's work to that of Blum, Kiefer and Rosenblatt does not appear to
have been noted before.

Replacing absolute values in $h$ by squares,
it is straightforward to show that a thus modified~$D$ reduces to
\begin{eqnarray*}
&&4 \bigl(E\bigl[F_1(X_1)-F_1(X_2)
\bigr] \bigl[F_2(Y_1)-F_2(Y_2)
\bigr] \bigr)^2
\\
&&\quad= 16 \bigl(E\bigl[F_1(X_1)-EF_1(X)
\bigr] \bigl[F_2(Y_1)-EF_2(Y)\bigr]
\bigr)^2 = \tfrac 19\tilde \rho_S^2 ,
\end{eqnarray*}
where
\[
\tilde\rho_S = 12 E\bigl[F_1(X_1)-EF_1(X)
\bigr] \bigl[F_2(Y_1)-EF_2(Y)\bigr]
\]
is a version of Spearman's correlation which coincides with $\rho_S$
given in Section~1 for continuous distributions (see Section 5 in
Kruskal \cite
{kruskal58}, for more details).

Following Kruskal's \cite{kruskal58} preference for Kendall's tau
over Spearman's rho due to its relative simplicity, the same preference
might be expressed for $\tau^*$ compared to $D$.

\subsection{Comparison to other ordinal consistent tests of independence}\label{sec-rel}

We now describe further approaches to obtaining consistent independence
tests for ordinal variables described in the literature.
It may be noted that $H$ and $D$ are special cases of a general family
of coefficients, which can be formulated as
%
\begin{equation}
\label{hoeff-gen} Q_{g,h} = Q_{g,h}(X,Y) = \int g
\bigl(\bigl|F_{12}(x,y)-F_1(x)F_2(y)\bigr|\bigr)\,\mathrm{d}
\bigl[h(F_{12}) (x,y)\bigr].
\end{equation}
For appropriately chosen $g$ and $h$, $Q_{g,h}=0$ if and only if $X$
and $Y$ are independent. Instances were studied by de~Wet \cite{dewet80},
Deheuvels \cite{deheuvels81}, Schweizer and Wolff \cite{sw81} and Feuerverger \cite{feuerverger93} (where the
former two focussed on asymptotic distributions of empirical versions,
while the latter two focussed on population coefficients).
Alternatively, R{\'e}nyi \cite{renyi59} proposed \emph{maximal correlation},
defined as
\[
\rho^+ = \sup_{g,h}\rho\bigl(g(X),h(Y)\bigr),
\]
where the supremum is taken over square integrable functions.
Though applicable to ordinal random variables, $\rho^+$ does not
utilize the ordinal nature of the variables. Furthermore, it is hard to
estimate, and has the drawback that it may equal one for distributions
arbitrarily ``close'' to independence (Kimeldorf and Sampson \cite{ks78}). An ordinal variant,
proposed by Kimeldorf and Sampson \cite{ks78}, was to maximize the correlation over
non-decreasing square integrable functions.

\subsection{Comparison to non-ordinal consistent tests of independence}\label{non-ord}

Recently Sz{\'e}kely, Rizzo and Bakirov \cite{srb07} introduced a consistent test of independence for
Euclidean random variables. With $\psi_{XY}$ the characteristic
function of the distribution of $(X,Y)\in\mR^p\times\mR^q$, and
$\psi
_X$ and $\psi_Y$ the characteristic functions of the corresponding
marginal distributions, they defined
%
\begin{equation}
\label{char} \dcov^2(X,Y) = \frac{1}{c_pc_q}\int
_{\mR^p\times\mR^q}\frac
{|\psi
_{XY}(s,t)-\psi_X(s)\psi_Y(t)|^2}{\Vert t\Vert ^{1+p}\Vert s\Vert ^{1+q}}\,\mathrm{d}s\,\mathrm{d}t,
\end{equation}
where $c_p$ and $c_q$ are constants.
It holds true that $\dcov^2(X,Y)\ge0$ with equality if and only if $X$
and $Y$ are independent, which is easy to show from the definition.
The expression~(\ref{char}) was originally introduced by Feuerverger \cite
{feuerverger93}, but only for real $X$ and $Y$ ($p=q=1$).

It was shown that $\dcov$ can equivalently be defined as
\[
\dcov^2(X,Y) = E\Vert X_1-X_2\Vert
\Vert Y_1-Y_2\Vert  + E\Vert X_1-X_2\Vert
E\Vert Y_1-Y_2\Vert  - 2E\Vert X_1-X_2\Vert
\Vert Y_1-Y_3\Vert.
\]
From this, it is straightforward to derive that
\[
\dcov^2(X,Y) = \tfrac1{4}Eh(X_1,X_2,X_3,X_4)h(Y_1,Y_2,Y_3,Y_4),
\]
where $h$ is defined by~(\ref{hdef}). Hence, for the case that $X$ and
$Y$ are real (i.e., $p=q=1$), $\dcov$ is closely related to $\tau^*$,
$\tau^*$ being a sign version.

With $Z_1$ and $Z_2$ independent with distribution $F$, let
\[
h_F(z_1,z_2) = -\tfrac12Eh(z_1,z_2,Z_1,Z_2),
\]
where $h$ is defined by~(\ref{hdef}). It can be verified that
\[
\dcov^2(X,Y) = Eh_{F_1}(X_1,X_2)h_{F_2}(Y_1,Y_2).
\]
As shown by Bergsma \cite{bergsma06} for the case that $X$ and $Y$ are real
and Sejdinovic \textit{et~al.} \cite{sgsf12} (see also Sejdinovic \textit{et~al.} \cite{ssgf12}) for the case that $X$ and
$Y$ are Euclidean, $h_F$ is a positive definite kernel implying
non-negativity of $\dcov^2$, while further properties of $h_F$ imply
equality to zero if and only if $X$ and $Y$ are independent. In fact,
as shown explicitly by Sejdinovic \textit{et~al.}, $\dcov^2$ falls in a
general class of association measures based on positive definite
kernels described by Gretton \textit{et~al.} \cite{gbss05}, which they called the
Hilbert--Schmidt independence criterion (HSIC). This criterion is a
generalization of Escoufier's vector covariance (Escoufier \cite
{escoufier73}, Robert and Escoufier \cite{re76}). It appears that $\tau^*$ is not an HSIC.

Although $\dcov^2$ and $\tau^*$ are similar in form, proofs of their
basic properties are very different. In particular, in spite of its
simple mathematical description, the proof for $\tau^*$ is much more
complex. The reason for this is that it appears hard to formulate $\tau
^*$ in terms of positive definite kernels, or as the expectation of a
squared norm of a random quantity (see also Lyons \cite{lyons12}).

Finally, another recent consistent test of independence for Euclidean
random variables is given by~Heller, Gorfine and Heller \cite{hhg12}, which is based on the
summation of Pearson chi-square statistics for well-chosen collapsing
of the bivariate distribution onto $2\times2$ contingency tables.

\section{Testing independence}\label{sec-applic}

A suitable test for independence is a permutation test which rejects
the independence hypothesis for large values of $t^*$, the empirical
value of $\tau^*$.
As an exact permutation test is too time consuming for moderately large
$n$, we use a Monte Carlo approximation, which is also called a
resampling test, and which is carried out as follows.
For $r=1,2,\ldots,$ let $(i_{r1},\ldots,i_{rn})$ be a random
permutation of $(1,\ldots,n)$, and let $t^*_r$ be $t^*$ computed for
the $r$th resample $(X_1,Y_{i_{r1}}),\ldots,(X_n,Y_{i_{rn}})$.
Then the Monte Carlo permutation $p$-value based on $R$ resamples is
computed as
\[
\mbox{Monte Carlo $p$-value} = \frac1R\sum_{r=1}^R
I\bigl(t^*_r>t^*\bigr).
\]
A further computational problem is the evaluation of $t^*$ itself (and
of the $t^*_r$), which requires computational time $\mathrm{O}(n^4)$, and may be
practically infeasible for moderately large samples. However, $t^*$ can
be well-approximated by taking a sufficiently large random sample of
subsets of four observations to approximate the sum in~(\ref{tdef2}).


As is well known, the permutation test conditions on the empirical
marginal distributions, which are sufficient statistics for the
independence model. In categorical data analysis, it is usually
referred to as an exact conditional test.
Note that there does not seem to be a need for an asymptotic
approximation to the sampling distribution of $t^*$.

\begin{table}[b]
\tablewidth=300pt
\caption{Artificial contingency table containing multinomial counts.
Permutation tests based on Kendall's tau and the Pearson chi-square
statistic do not yield a significant association ($p=0.99$, resp.,
$p=0.25$), but a permutation test based on $t^*$ yields $p=0.035$}\label{ct1}
\begin{tabular*}{300pt}{@{\extracolsep{\fill}}llllllll@{}}
\hline
&\multicolumn{7}{l@{}}{$Y$}\\[-6pt]
&\multicolumn{7}{l@{}}{\hrulefill}\\
 $X$& 1 & 2 & 3 & 4 & 5 & 6 & 7 \\
\hline
1 & 2 & 1 & 0 & 0 & 0 & 1 & 2 \\
2 & 1 & 2 & 0 & 0 & 0 & 2 & 1 \\
3 & 0 & 0 & 2 & 1 & 2 & 0 & 0 \\
4 & 0 & 0 & 1 & 1 & 1 & 0 & 0 \\
5 & 0 & 0 & 1 & 2 & 1 & 0 & 0 \\
\hline
\end{tabular*}
\end{table}

In this section, we compare various tests of independence using an
artificial and a real data set and via a simulation study.

\subsection{Examples}

An artificial multinomial table of counts is given in Table~\ref{ct1},
where $X$ and $Y$ are ordinal variables with~5 and~7 categories.
Visually, we can detect an association pattern, but as it is
non-monotonic a test based on Kendall's tau does not yield a
significant $p$-value. The chi-square test also yields a
non-significant $p=0.252$, while a permutation test based on $t^*$
yields $p=0.032$, giving evidence of an association. We also did tests
based on $D$, which yields $p=0.047$, and the test based on Hoeffding's
$H$ yields $p=0.028$.
In this example, using a consistent test designed for ordinal data,
evidence for an association can be found, which is not possible with a
nominal data test like the chi-square test or with a test based on
Kendall's tau. For all tests except Hoeffding's $R=10^6$ resamples were
used, and for Hoeffding's test $R=4000$ resamples were used.

\begin{table}
\tablewidth=300pt
\caption{Results of study comparing two treatments of gastric ulcer}
\label{ct2}
\begin{tabular*}{300pt}{@{\extracolsep{\fill}}lllll@{}}
\hline
& \multicolumn{4}{l}{Change in size of Ulcer Crater ($Y$)}\\[-6pt]
& \multicolumn{4}{l@{}}{\hrulefill}\\
Treatment group ($X$) &  Larger & \multicolumn{1}{c}{Healed ($<\frac23$)} & \multicolumn{1}{c}{Healed ($\ge\frac23$)} & Healed \\
\hline
 $A$ & \phantom{0}6 & 4 & 10 & 12 \\
$B$ & 11 & 8 & \phantom{0}8 & \phantom{0}5 \\
\hline
\end{tabular*}
\end{table}

Table~\ref{ct2} shows data from a randomized study to compare two
treatments for a gastric ulcer crater, and was previously analyzed in
Agresti \cite{agresti10}. Using $R=10^5$ resamples, the chi-square test yields
$p=0.118$, Kendall's tau yields $p=0.019$, $t^*$ yields $p=0.028$, $D$
yields $p=0.026$, and using $10^4$ resamples Hoeffding's $H$ yields $p=0.006$.

\subsection{\texorpdfstring{Simulated average $p$-values for independence tests based on $D$, $H$, and~$\tau^*$}
{Simulated average $p$-values for independence tests based on $D$, $H$, and tau*}}

Any of the three tests can be expected to have most power of the three
for certain alternatives, and least power of the three for others.
Given the broadness of possible alternatives, it cannot be hoped to get
a simple description of alternatives for which any single test is the
most powerful. However, some insight may be gained by looking at
average $p$-values for a set of carefully selected alternatives.

\begin{figure}

\includegraphics{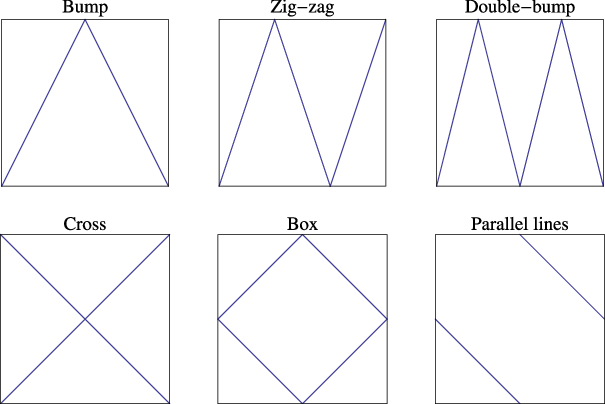}

\caption{Simulations were done for data generated from the uniform
distribution on the lines within each of the six boxes. For all except
the Zig--zag and the Parallel lines, the ordinary correlation is
zero.}\label{simpats}
\end{figure}

In Figure~\ref{simpats}, six boxes with lines in them are represented,
and we simulated from the uniform distribution on these lines. The
first five maximize or minimize the correlation between some simple
orthogonal functions for given uniform marginals. In particular, say
the boxes represent the square $[0,1]\times[0,1]$, then the Bump,
Zig--zag and Double bump distributions maximize, for given uniform marginals,
\[
\rho\bigl[\cos(2\pi X),\cos(\pi Y)\bigr],\qquad \rho\bigl[\cos(3\pi X),\cos(\pi Y)
\bigr]\quad\mbox{and}\quad \rho\bigl[\cos(4\pi X),\cos(\pi Y)\bigr],
\]
respectively. The Cross and Box distributions respectively maximize and
minimize, for given uniform marginals,
\[
\rho\bigl[\cos(2\pi X),\cos(2\pi Y)\bigr].
\]
As they represent in this sense extreme forms of association, these
distributions should yield good insight in the comparative performance
of the tests.
Furthermore, the Parallel lines distribution was chosen because it is
simple and demonstrates a weakness of Hoeffding's test, as it has
comparatively very little power here (we did not manage to find a
distribution where $D$ or $\tau^*$ fare so comparatively poorly).
Note that all six distributions have uniform marginals and so are
copulas, and several were also discussed in Nelsen~\cite{nelsen06}.

\begin{figure}

\includegraphics{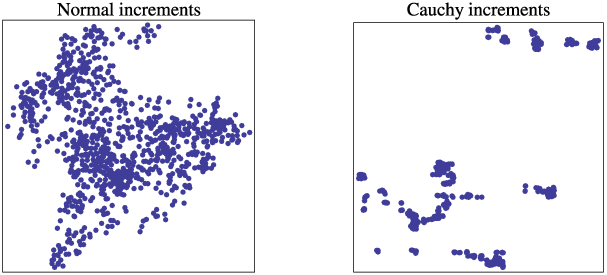}

\caption{1000 points of a random walk. In the first plot the $(x,y)$
increments are independent normals, in the second they are independent
Cauchy variables.}\label{simpats2}
\end{figure}

We also did a Bayesian simulation, based on random distributions with
dependence. In particular, the data are $(X_1,Y_1),\ldots,(X_n,Y_n)$,
where, for i.i.d. $(\varepsilon_{1i},\varepsilon_{2i})$,
\begin{eqnarray*}
(X_1,Y_1) &=& (\varepsilon_{11},
\varepsilon_{21}),
\\
(X_{i+1},Y_{i+1})&=& (X_{i},Y_{i})+(
\varepsilon_{1i},\varepsilon_{2i})\qquad i=1,\ldots,n-1.
\end{eqnarray*}
Of course, the $(X_i,Y_i)$ are not i.i.d., but conditioning on the
empirical marginals the permutations of the $Y$-values give equally
likely data sets under the null hypothesis of independence, so the
permutation test is valid.
Two distributions for the increments $(\varepsilon_{1i},\varepsilon
_{2i})$ were used: independent normals and independent Cauchy distributions.
In Figure~\ref{simpats2}, points generated in this way are plotted.
Note that for the Cauchy increments, the heavy tails of the marginal
distributions are automatically taken care of by the use of ranks, so
in that respect the three tests described here are particularly suitable.

Finally, we also simulated normally distributed data with correlation~$0.5$.

\begin{table}[b]
\caption{Average $p$-values. See Figures~\protect\ref{simpats} and
\protect\ref{simpats2} and the text for explanations}\label{simstudy1}
\begin{tabular*}{\textwidth}{@{\extracolsep{\fill}}lllll@{}}
\hline
 &  & \multicolumn{3}{l}{Average
$p$-value}\\[-6pt]
&&\multicolumn{3}{l@{}}{\hrulefill}\\
Distribution& Sample size $n$& $D$ & $H$ & $\tau^*$ \\
\hline
Random walk (normal increments) & 50 & 0.061 & 0.080 & 0.061 \\
Random walk (Cauchy increments) & 30 & 0.039 & 0.065 & 0.031 \\[3pt]
Bump & 12 & 0.087 & 0.061 & 0.045 \\
Zig--zag & 25 & 0.083 & 0.011 & 0.036 \\
Double-bump & 30 & 0.056 & 0.005 & 0.019 \\
Cross & 50 & 0.052 & 0.003 & 0.021 \\
Box & 50 & 0.070 & 0.008 & 0.019 \\
Parallel lines & 10 & 0.055 & 0.710 & 0.076 \\[3pt]
Normal distribution ($\rho=0.5$) & 30 & 0.055 & 0.052 & 0.073 \\
\hline
\end{tabular*}
\end{table}

Average $p$-values are given in Table~\ref{simstudy1}, where all
averages are over at least 40,000 simulations (for $D$, we did 200,000
simulations). Hoeffding's test compares extremely badly with our test
for the parallel lines distribution, and is worse than our test for the
random walks, but outperforms our test for the Zig--zag, Double-bump,
Cross and Box distributions.
The reason for the poor performance of Hoeffding's test for the
parallel lines distribution is that five points can only be concordant
(see Section~\ref{sec-hoeff}) if they all lie on a single line (a
discordant set of five points has zero probability). Similarly, for the
Zig--zag, Double-bump and Cross concordant sets of five points can be
seen to be especially likely, so these choices of distributions favour
the Hoeffding test.
Note that Hoeffding's test is less suitable for general use because it
is not necessarily zero under independence if there is a positive
probability of tied observations.

The BKR test fares slightly worse than ours for the random walk with
Cauchy increments, and significantly worse than ours for the Bump,
Zig--zag, Cross and Box distributions, and does somewhat better than
ours for the normal distribution.
It appears that the BKR test has more power than ours for a monotone
alternative (such as the normal distribution), at the cost of less
power for some more complex alternatives.

\section{\texorpdfstring{Proof of Theorem~\protect\ref{th1}}
{Proof of Theorem 1}}\label{sec3}
Here we give the proof of Theorem~\ref{th1} for arbitrary real random
variables $X$ and $Y$. A shorter proof for continuous $X$ and $Y$ is
given by Dassios and Bergsma \cite{db12}. Readers wishing to gain an understanding of the
essence of the proof may wish to study the shorter proof first.

First, consider three real valued random variables $U$, $V$and $W$.
They have
continuous densities $\tilde{f}  ( x  )$, $\tilde{g}
( x
)$ and \ $\tilde{k}  ( x  )$ as well as probability
masses $f
( x_i  )$, $g  ( x_i  )$ and $k  ( x_i
)$ at
points $x_1, x_2, \ldots.$ We also define
\begin{eqnarray*}
F ( x )& =& P ( U < x ) = \sum_{x_i < x} f (
x_i ) + \int_{y < x} \tilde{f} ( y ) \,\tmop{d}y,
\\
G ( x ) &=& P ( V < x ) = \sum_{x_i < x} g (
x_i ) + \int_{y < x} \tilde{g} ( y ) \,\tmop{d}y
\end{eqnarray*}
and
\[
K ( x ) = P ( W < x ) = \sum_{x_i < x} k (
x_i ) + \int_{y < x} \tilde{k} ( y ) \,\tmop{d}y .
\]
We will also use $H  ( x  ) = \frac{K  ( x  )}{G
( x
)}$. Note that $H  ( x  )$ also admits the representation
\[
H ( x ) = \sum_{x_i < x} h ( x_i ) + \int
_{y < x} \tilde{h} ( y ) \,\tmop{d}y
\]
but unlike the other three function that are non-decreasing $\tilde{h}
(
x  )$ and $h  ( x_i  )$ can take negative values.

We start by proving the following intermediate result.

\begin{lemma}
\label{lem-ins3}Assume that $G  ( x  ) = 1$ implies $F
( x
) =K  ( x
) = 1$ and that there is a constant $c$ such that $F
( x  ) \leq c G  ( x  )$ and $K  ( x
) \leq c G^2  ( x  )$ for all $x$. Define
\begin{eqnarray*}
S &=& 2 \sum \bigl( F ( x_i ) - G ( x_i )
\bigr) \bigl( F ( x_i ) g ( x_i ) - G ( x_i
) f ( x_i ) \bigr) \frac{K  ( x_i  )}{G^2
( x_i
)}\\
&&{} -
\sum \bigl( F ( x_i ) g ( x_i ) - G (
x_i ) f ( x_i ) \bigr)^2 \frac{K  ( x_i  )}{G^2
( x_i  )}\\
&&{}+
2 \int_{\nosymbol} \bigl( F ( x ) - G ( x ) \bigr) \bigl( F ( x )
\tilde{g} ( x ) - G ( x ) \tilde{f} ( x ) \bigr) \frac{K  (
x_{\nosymbol}
)}{G^2  ( x_{\nosymbol}  )} \,\tmop{d}x,
\end{eqnarray*}
where summation is over all $x_i$ such that $K  ( x_i  ) > 0$
and at least one of $f  ( x_i  )$ and $g  ( x_i
)$ is
positive, and
integration is over all $x$ such that $K  ( x
) > 0$.

We then have $S \geq0$ with equality iff $F  ( x  ) = G
( x
)$ for all $x_{\nosymbol}$ such that $K  ( x_{\nosymbol}
) >
0$.
\end{lemma}
\begin{pf}
The conditions stated in the lemma ensure that the sums and
integral exist. We can rewrite
\begin{eqnarray*}
S &=& 2 \sum \bigl( F ( x_i ) - G ( x_i )
\bigr) \bigl( F ( x_i ) g ( x_i ) - G ( x_i
) f ( x_i ) \bigr) \frac{H  ( x_i  )}{G  ( x_i
)}\\
&&{} -
\sum \bigl( F ( x_i ) g ( x_i ) - G (
x_i ) f ( x_i ) \bigr)^2 \frac{H  ( x_i  )}{G
( x_i  )}
\\
&&{}+
2 \int_{\nosymbol} \bigl( F ( x ) - G ( x ) \bigr) \bigl( F ( x )
\tilde{g} ( x ) - G ( x ) \tilde{f} ( x ) \bigr) \frac{H  ( x_{\nosymbol}
)}{G
( x_{\nosymbol}  )} \,\tmop{d}x .
\end{eqnarray*}
For simplicity, we denote $F  ( x  ), G  ( x  ), H
( x
), f  ( x_i  ), g  ( x_i  ), h  ( x_i
),
\tilde{f}  ( x  ), \tilde{g}  ( x  )$ and
$\tilde{h}
(
x  )$ by $F, G,  H, f, g, h, \tilde{f}, \tilde{g}$ and $\tilde
{h}$. We
have
\begin{eqnarray}\label{11}
S& =& 2 \sum_{\nosymbol} ( F - G ) \bigl( ( F - G ) g - G (
f - g ) \bigr) \frac{H}{G} \nonumber\\
&&{}+
2 \int ( F - G ) \bigl( ( F - G ) \tilde{g} - G ( \tilde{f} - \tilde{g} ) \bigr)
\frac{H}{G} \,\tmop {d}x\nonumber\\
&&{} -
\sum_{\nosymbol} \bigl( ( F - G ) g - G ( f - g )
\bigr)^2 \frac{H}{G}
\nonumber\\
& =&
2 \sum ( F - G )^2 \frac{H}{G} g + 2 \int
_{\nosymbol} ( F - G )^2 \frac{H}{G} \tilde{g}
\,\tmop{d}x\\
&&{} -
2 \sum_{\nosymbol} H ( F - G )^{\nosymbol} ( f - g )\nonumber\\
&&{} -
2 \int_{\nosymbol} H ( F - G )^{\nosymbol} ( \tilde {f} -
\tilde{g} ) \,\tmop{d}x\nonumber\\
&&{} -
\sum \bigl( ( F - G ) g - G ( f - g ) \bigr)^2
\frac{H}{G}^{\nosymbol}.\nonumber
\end{eqnarray}
The function $H  ( F - G  )^2$ vanishes at $- \infty$ (because
of the
conditions of the lemma) and $+ \infty$. Considering its integral and sum
representation we have
\begin{eqnarray*}
&& 2 \sum H ( F - G )^{\nosymbol} ( f - g ) + 2 \int
_{\nosymbol} H ( F - G )^{\nosymbol} ( \tilde {f} - \tilde{g} )
\,\tmop{d}x\\
&&\quad{} +
\sum ( F - G )^2 h + \int_{\nosymbol} ( F - G
)^2 \tilde{h} \,\tmop{d}x \\
&&\quad{}+ 2 \sum ( F - G )^{\nosymbol} ( f - g ) h + \sum
^{\nosymbol} ( f - g )^2 h + \sum
^{\nosymbol} H ( f - g )^2 = 0 ,
\end{eqnarray*}
and therefore
\begin{eqnarray}\label{2}
&&- 2 \sum H ( F - G )^{\nosymbol} ( f - g ) - 2 \int
_{\nosymbol} H ( F - G )^{\nosymbol} ( \tilde {f} - \tilde{g} )
\,\tmop{d}x\nonumber \\
&&\quad =
\sum ( F - G )^2 h + \int_{\nosymbol} ( F - G
)^2 \tilde{h} \,\tmop{d}x\\
&&\qquad{}
+ 2 \sum ( F - G )^{\nosymbol} ( f - g ) h + \sum
^{\nosymbol} ( f - g )^2 h + \sum
^{\nosymbol} H ( f - g )^{2 \nosymbol} .\nonumber
\end{eqnarray}
Moreover,
\begin{eqnarray}\label{3}
&&\frac{H}{G} \bigl( ( F - G ) g - G ( f - g ) \bigr)^2
\nonumber
\\[-8pt]
\\[-8pt]
\nonumber
&&\quad=
( F - G )^2 g^2 \frac{H}{G} + G^{\nosymbol} H
( f - g )^2 - 2 ( F - G ) ( f - g ) Hg.
\end{eqnarray}
Substituting (\ref{2}) \ and (\ref{3}) into (\ref{11}), and denoting $M
= F -
G$, $m = f - g$ and $\tilde{m} = \tilde{f} - \tilde{g}$ we have
\begin{eqnarray*}
S& =& \sum M^2 \biggl( 2 g \frac{H}{G} + h -
g^2 \frac{H}{G} \biggr) + 2 \sum Mm ( h + gH ) + \sum
m^2 ( H + h - GH ) \\
&&{}+
\int M^2 \biggl( 2 \tilde{g} \frac{H}{G} + \tilde{h} \biggr)
\,\tmop {d}x\\
& =&
\sum ( M + m )^2 \biggl( g \frac{H}{G + g} + h \biggr) +
\sum M^2 \biggl( 2 g \frac{H}{G} - g
\frac{H}{G + g} - g^2 \frac{H}{G} \biggr)\\
&&{} -
2 \sum Mm \biggl( g \frac{H}{G + g} - gH \biggr) + \sum
m^2 \biggl( H - GH - g \frac{H}{G + g} \biggr)\\
&&{} +
\int M^2 \biggl( \tilde{g} \frac{H}{G} + \tilde{h} \biggr)
\,\tmop {d}x + \int M^2 \tilde{g} \frac{H}{G} \,\tmop{d}x\\
& =&
\sum ( M + m )^2 \biggl( g \frac{H}{G + g} + h \biggr) +
\int M^2 \biggl( \tilde{g} \frac{H}{G} + \tilde{h} \biggr)
\,\tmop{d}x + \int M^2 \tilde{g} \frac{H}{G} \,\tmop{d}x\\
&&{} +
\sum M^2 \biggl( g \frac{H}{G} + g^2
\frac{H  ( 1 - G - g  )}{G
( G + g  )} \biggr) - 2 \sum Mm \biggl( g \frac{H  (
1 -
G - g
)}{G + g}
\biggr)\\
&&{} +
\sum m^2 \frac{H}{G + g} \bigl( ( 1 - G ) G - g G \bigr)
.
\end{eqnarray*}
Observe now that since $K = HG$
\[
g \frac{H}{G + g} + h = \frac{gH + hG + hg}{G + g} = \frac{k}{G + g} \geq0
\]
and
\[
\tilde{g} \frac{H}{G} + \tilde{h} = \frac{\tilde{k}}{G} \geq0.
\]
Moreover, the quadratic form
\begin{eqnarray*}
&& M^2 \biggl( g \frac{H}{G} + g^2
\frac{H  ( 1 - G - g  )}{G
(
G + g  )} \biggr) - 2 Mm \biggl( g \frac{H  ( 1 - G - g
)}{G
+ g} \biggr) +
m^2 \frac{H}{G + g} \bigl( ( 1 - G ) G - gG \bigr)\\
&&\quad =
\frac{M^2 g H}{G} + ( M g - m G )^2 \frac{H  ( 1 -
G - g
)}{G  ( G + g  )}.
\end{eqnarray*}

All terms in $S$ are non-negative and are equal to zero iff $M  ( x
) = 0$ for all $x$ such that $K  ( x
) > 0 $, that
is the two distributions $F$ and $G$ are identical for all
$x_{\nosymbol
}$ such
that $K  ( x_{\nosymbol}  ) > 0$.
\end{pf}

Before we prove Theorem \ref{th1}, we will prove another result as it will be used
repeatedly.

\begin{lemma}
\label{lem-ins4}Let $A$, $B$ and $C$ be events in the same probability space
as the random variable $X$ and define
\begin{eqnarray*}
L \bigl( x^{ ( 1  )}, x^{ ( 2  )} \bigr) &=& \bigl( P \bigl( A \vert X
= x^{ ( 1  )}  \bigr) - P \bigl( A \llvert X < x^{ ( 1  )}
\wedge x^{ ( 2  )} \bigr) \bigr) \\
&&{}\times
\bigl( P \bigl( A \vert X = x^{ ( 2  )}  \bigr) - P \bigl( A
\llvert X < x^{ ( 1  )} \wedge x^{ ( 2  )}  \bigr) \bigr)\\
&&{} \times
P \bigl( B \vert X < x^{ ( 1  )} \wedge x^{ ( 2
)}  \bigr)
P \bigl( C \vert X < x^{ ( 1  )} \wedge x^{ ( 2  )}  \bigr)
\bigl( P \bigl( X < x^{
( 1
)} \wedge x^{ ( 2  )} \bigr)
\bigr)^2 .
\end{eqnarray*}
We then have
\[
E \bigl( L ( X_1, X_2 ) \bigr) \geq0
\]
with equality iff $P  ( X < x  ) = P  ( X < x \vert A
)$ \ for all $x$ such that \mbox{$P  ( X < x \vert  B_{\nosymbol}  ) P  ( X < x_{\nosymbol} \vert  C
) > 0$}.
\end{lemma}
\begin{pf}
Let $X$ have continuous density $\tilde{g}  ( x  )$ and
probability masses $g  ( x_i  )$ at points $x_1, x_2, \ldots$
and let
$X$ have continuous density $\tilde{g}_A  ( x  )$ and probability
masses $g_A  ( x_i  )$ at points $x_1, x_2, \ldots$
conditionally on
$A$. Define also
\[
G ( x ) = P ( X < x ) = \sum_{x_i < x} g (
x_i ) + \int_{y < x} \tilde{g} ( y ) \,\tmop{d}y
\]
and
\[
G_A ( x ) = P \bigl( X < x \vert  A \bigr) = \sum
_{x_i
< x} g_A ( x_i ) + \int
_{y < x} \tilde{g}_A ( y ) \,\tmop{d}y .
\]
Conditioning on values of $X_1 \wedge X_2$ and using Bayes' theorem, we can
see that
\begin{eqnarray*}
&& E \bigl( L ( X_1, X_2 ) \bigr) \\
&&\quad= \bigl( P ( A )
\bigr)^2 \sum P \bigl( B \vert X < x_i
 \bigr) P \bigl( C \vert X < x_i  \bigr)\\
 &&\qquad{} \times
\bigl\{ 2 \bigl( \bigl( 1 - G_A ( x_i ) \bigr) G (
x_i ) - \bigl( 1 - G ( x_i ) \bigr) G_A (
x_i ) \bigr) \bigl( g_A ( x_i ) G (
x_i ) - g ( x_i ) G_A ( x_i )
\bigr)\\
&&\hspace*{15pt}\qquad{} -
\bigl( g_A ( x_i ) G ( x_i ) - g (
x_i ) G_A ( x_i ) \bigr)^2 \bigr\}
\\[-2pt]
&&\qquad{} +
\bigl( P ( A ) \bigr)^2 \int P \bigl( B \vert X < x_{\nosymbol}
 \bigr) P \bigl( C \vert X < x_{\nosymbol}  \bigr)\\[-2pt]
 &&\hspace*{55pt}\qquad{} \times 2 \bigl( \bigl( 1 - G_A ( x_{\nosymbol} ) \bigr) G (
x_{\nosymbol} ) - \bigl( 1 - G ( x_{\nosymbol} ) \bigr) G_A (
x_{\nosymbol} ) \bigr) \bigl( \tilde{g}_A ( x ) G ( x ) -
\tilde{g} ( x ) G_A ( x_{\nosymbol} ) \bigr) \,\tmop{d}x\\[-2pt]
&&\quad =
P ( B ) P ( C ) \bigl( P ( A ) \bigr)^2 \sum
\frac{K  ( x_i  )}{G^2  ( x_i  )} \cdot
\bigl\{ 2 \bigl( G ( x_i ) - G_A ( x_i )
\bigr) \bigl( g_A ( x_i ) G ( x_i ) - g (
x_i ) G_A ( x_i ) \bigr)\\[-2pt]
&&\hspace*{138pt}\qquad{}  -  \bigl(
g_A ( x_i ) G ( x_i ) - g (
x_i ) G_A ( x_i ) \bigr)^2
\bigr\} \\[-2pt]
&&\qquad{} +
P ( B ) P ( C ) \bigl( P ( A ) \bigr)^2 \int\frac{K  ( x_{\nosymbol}  )}{G^2  ( x_{\nosymbol}
)} 2
\bigl( G ( x_{\nosymbol} ) - G_A ( x_{\nosymbol} ) \bigr)
\bigl( \tilde{g}_A ( x ) G ( x ) - \tilde{g} ( x ) G_A (
x_{\nosymbol} ) \bigr) \,\tmop{d}x,
\end{eqnarray*}
where
\[
K ( x_{\nosymbol} ) = P \bigl( X < x \vert B _{\nosymbol} \bigr)
P \bigl( X < x_{\nosymbol} \vert C  \bigr) \nosymbol.
\]
The result then follows from Lemma~\ref{lem-ins3} ($F = G_A$). It is
easy to
see that the conditions in Lemma~\ref{lem-ins3} are satisfied. For example, $P
( X < x
\vert  B   ) \leq\frac{P  ( X < x  )}{P
( B  )}$.
\end{pf}


\begin{pf*}{Proof of Theorem \ref{th1}}
We need to prove that
\begin{eqnarray*}
&& P ( Y_1 \wedge Y_2 > Y_3 \vee
Y_4, X_3 \vee X_4 < X_1 \wedge
X_2 ) \\[-2pt]
&&\quad{}+
P ( Y_1 \vee Y_2 < Y_3 \wedge
Y_4, X_3 \vee X_4 < X_1 \wedge
X_2 )\\[-2pt]
&&\quad{} -
P ( Y_1 \wedge Y_3 > Y_2 \vee
Y_4, X_3 \vee X_4 < X_1 \wedge
X_2 )\\[-2pt]
&&\quad{} -
P ( Y_1 \vee Y_3 < Y_2 \wedge
Y_4, X_3 \vee X_4 < X_1 \wedge
X_2 ) \geq0
\end{eqnarray*}
with equality in the independence case.

Let $ ( X, Y  )$ represent any of the pairs $ ( X_{i,} Y_i
)$. Define now $F_1  ( y  ) = P  ( Y_{\nosymbol
} < y
\vert
X_{\nosymbol} = x_{\nosymbol}^{ ( 1  )}
)$, \ $F_2
( y  ) = P  ( Y_{\nosymbol} < y \vert  X_{\nosymbol
} =
x_{\nosymbol}^{ ( 2  )}   )$ and $G  ( y
) = P
( Y_{\nosymbol} < y \vert  X_{\nosymbol} < x_{\nosymbol
}^{ ( 1
)} \wedge x_{\nosymbol}^{ ( 2  )}  )$
\
with the
representations
\begin{eqnarray*}
F_1 ( y )& =& \sum_{y_i < y} f_1
( y_i ) + \int_{z
< y} \widetilde{f}_1
( z ) \,\tmop{d}z,
\\[-2pt]
F_2 ( y )& =& \sum_{y_i < y} f_2 (
y_i ) + \int_{z
< y} \widetilde{f}_2 (
z ) \,\tmop{d}z
\end{eqnarray*}
and
\[
G ( y ) = \sum_{y_i < y} g ( y_i ) + \int
_{z < y} \tilde{g} ( z ) \,\tmop{d}z .
\]
Note that conditionally on the event
\[
\Theta= \bigl\{ X_1 = x_{\nosymbol}^{ ( 1  )},
X_2 =_{\nosymbol} x_{\nosymbol
}^{
( 2
)},
X_3 < x_{\nosymbol}^{ ( 1  )} \wedge x_{\nosymbol
}^{ (
2  )},
X_4 < x_{\nosymbol}^{ ( 1  )} \wedge x_{\nosymbol}^{ ( 2  )}
\bigr\},
\]
the distribution of the minimum of
$Y_1$ and $Y_2$ has density $ ( 1 - F_1  ) \widetilde{f}_2 +
( 1
- F_2  ) \widetilde{f}_1$ and probability masses $ ( 1 - F_1
)
f_2 +  ( 1 - F_2  ) f_1 - f_1 f_2$ at $y_1, y_2, \ldots,$ the
distribution of the minimum of \ $Y_3$ and $Y_4$ has density\vadjust{\goodbreak} $2  (
1 - G
) \tilde{g}$ and probability masses $2  ( 1 - G  ) g
- g^2
\nosymbol$, the \ distribution of the minimum of $Y_1$ and $Y_3$ has density
$ ( 1 - F_1  ) \tilde{g} +  ( 1 - G  )
\widetilde
{f}_1$ and
probability masses $ ( 1 - F_1  ) g +  ( 1 - G  ) f_1
- f_1
g$ and the distribution of the minimum of $Y_2$ and $Y_4$ has density
$ (
1 - F_2  ) \tilde{g} +  ( 1 - G  ) \widetilde{f}_2$ and
probability masses $ ( 1 - F_2  ) g +  ( 1 - G  ) f_2
- f_2
g$. We therefore have (suppressing the arguments of the functions)
\begin{eqnarray*}
&& P \bigl( Y_1 \wedge Y_2 > Y_3 \vee
Y_4 \vert \Theta \bigr) + P \bigl( Y_1 \vee
Y_2 < Y_3 \wedge Y_4 \vert \Theta
\bigr)\\[-2pt]
&&\qquad{} -
P \bigl( Y_1 \wedge Y_3 > Y_2 \vee
Y_4 \vert \Theta \bigr) - P \bigl( Y_1 \vee
Y_3 < Y_2 \wedge Y_4 \llvert \Theta
\bigr) \\[-2pt]
&&\quad=
\sum \bigl( ( 1 - F_1 ) f_2 + ( 1 -
F_2 ) f_1 - f_1 f_2 \bigr)
G^2 + \sum \bigl( 2 ( 1 - G ) g - g^2 \bigr)
F_1 F_2 \\[-2pt]
&&\qquad{}-
\sum \bigl( ( 1 - F_1 ) g + ( 1 - G ) f_1 -
f_1 g \bigr) F_2 G - \sum \bigl( ( 1 -
F_2 ) g + ( 1 - G ) f_2 - f_2 g \bigr)
F_1 G\\[-2pt]
&&\qquad{} +
\int \bigl( ( 1 - F_1 ) \widetilde{f}_2 + ( 1 -
F_2 ) \widetilde{f}_1 \bigr) G^2 \,\tmop{d}y +
\int2 ( 1 - G ) \tilde g F_1 F_2 \,\tmop{d}y\\[-2pt]
&&\qquad {}-
\int \bigl( ( 1 - F_1 ) \tilde{g} + ( 1 - G ) \widetilde{f}_1
\bigr) F_2 G\,\tmop{d}y - \int \bigl( ( 1 - F_2 ) \tilde{g}
+ ( 1 - G ) \widetilde{f}_2 \bigr) F_1 G\,\tmop{d}y \\[-2pt]
&&\quad=
\sum ( F_1 - G ) ( F_2 g - Gf_2 )
+ \sum ( F_2 - G ) ( F_1 g - Gf_1
) - \sum ( F_1 g - Gf_1 ) ( F_2 g
- Gf_2 )\\[-2pt]
&&\qquad{} +
\int ( F_1 - G ) ( F_2 \tilde{g} - G
\widetilde{f}_2 ) \,\tmop{d}y + \int ( F_2 - G ) (
F_1 \tilde {g} - G \tilde{f}_1 ) \,\tmop{d}y\\[-2pt]
&&\quad =
2 \sum ( F_1 - G ) ( F_2 - G ) g - \sum
( F_1 - G ) ( f_2 - g ) G - \sum (
F_2 - G ) ( f_1 - g ) G\\[-2pt]
&&\qquad{} -
\sum ( F_1 g - Gf_1 ) ( F_2 g -
Gf_2 ) + 2 \int ( F_1 - G ) ( F_2 - G )
\tilde{g} \,\tmop{d}y \\[-2pt]
&&\qquad{}-
\int ( F_1 - G ) ( \tilde{f}_2 - \tilde{g} ) G\,\tmop{d}y
- \int ( F_2 - G ) ( \tilde{f}_1 - \tilde{g} )G
\,\tmop{d}y.
\end{eqnarray*}
The function $G  ( F_1 - G  )^{\nosymbol}  ( F_2 - G
)$
vanishes at $- \infty$ \ and $+ \infty$. Considering its integral and sum
representation, we have
\begin{eqnarray}\label{bb}
&& - \sum ( F_1 - G ) ( f_2 - g ) G - \sum
( F_2 - G ) ( f_1 - g ) G\nonumber\\[-2pt]
&&\qquad{} -
\int ( F_1 - G ) ( \tilde{f}_2 - \tilde{g} ) G\,\tmop{d}y
- \int ( F_2 - G ) ( \tilde{f}_1 - \tilde{g} ) G
\,\tmop{d}y \nonumber\\[-2pt]
&&\quad=
\sum ( F_1 - G ) ( F_2 - G ) g + \sum
( F_1 - G ) ( f_2 - g ) g + \sum (
F_2 - G ) ( f_1 - g ) g
\nonumber
\\[-2pt]
&&\qquad{} +
\sum ( f_1 - g ) ( f_2 - g ) G + \sum
( f_2 - g ) ( f_1 - g ) g \\[-2pt]
&&\qquad{}+ \int ( F_1 - G
) ( F_2 - G ) \tilde{g} \,\tmop{d}y\nonumber\\[-2pt]
&&\quad =
\sum ( F_1 + f_1 - G - g ) ( F_2
+ f_2 - G - g ) g\nonumber\\[-2pt]
&&\qquad{} +
\sum ( f_1 - g ) ( f_2 - g ) G + \int (
F_1 - G ) ( F_2 - G ) \tilde{g} \,\tmop{d}y
.^{\nosymbol}\nonumber
\end{eqnarray}
Moreover,
\begin{eqnarray}\label{7}
( F_1 g - Gf_1 ) ( F_2 g -
Gf_2 ) &=& ( F_1 - G ) ( F_2 - G )
g^2 - ( F_1 - G ) ( f_2 - g ) Gg \nonumber\\
&&{}-
( F_2 - G ) ( f_1 - g ) Gg + ( f_1 - g ) (
f_2 - g ) G^2\nonumber \\
&=&
( F_1 - G ) ( F_2 - G ) g^2 + (
f_1 - g ) ( f_2 - g ) G^2\nonumber\\
&&{} - (
F_1 + f_1 - G - g ) ( F_2 + f_2
- G - g ) gG \\
&&{}+
( F_1 - G ) ( F_2 - G ) gG + ( f_1 - g ) (
f_2 - g ) gG\nonumber \\
&=&
( F_1 - G ) ( F_2 - G ) g ( G + g ) + (
f_1 - g ) ( f_2 - g ) G ( G + g )\nonumber\\
&&{} -
( F_1 + f_1 - G - g ) ( F_2 +
f_2 - G - g ) gG.^{\nosymbol}\nonumber
\end{eqnarray}
Using (\ref{7}) and (\ref{bb}), we have
\begin{eqnarray*}
&& P \bigl( Y_1 \wedge Y_2 > Y_3 \vee
Y_4 \vert \Theta \bigr) + P \bigl( Y_1 \vee
Y_2 < Y_3 \wedge Y_4 \vert \Theta
\bigr)\\
&&\qquad{} -
P \bigl( Y_1 \wedge Y_3 > Y_2 \vee
Y_4 \vert \Theta \bigr) - P \bigl( Y_1 \vee
Y_3 < Y_2 \wedge Y_4 \vert \Theta
\bigr)\\
&&\quad =
\sum ( F_1 - G ) ( F_2 - G ) g + \sum
( F_1 - G ) ( F_2 - G ) g ( 1 - G - g ) \\
&&\qquad{}+
\sum ( F_1 + f_1 - G - g ) ( F_2
+ f_2 - G - g ) g \\
&&\qquad{}+ \sum ( F_1 +
f_1 - G - g ) ( F_2 + f_2 - G - g ) gG \\
&&\qquad{}+
\sum ( f_1 - g ) ( f_2 - g ) G ( 1 - G - g ) +
3 \int ( F_1 - G ) ( F_2 - G ) \tilde{g} \,\tmop{d}y .
\end{eqnarray*}
We therefore conclude that conditionally on \ $ \{ X_1 =
x_{\nosymbol}^{ ( 1  )}, X_2 =_{\nosymbol} x_{\nosymbol
}^{
( 2
)}  \}$,
\begin{eqnarray*}
&&P ( Y_1 \wedge Y_2 > Y_3 \vee
Y_4, X_3 \vee X_4 < X_1 \wedge
X_2 )\\
&&\qquad{} +
P ( Y_1 \vee Y_2 < Y_3 \wedge
Y_4, X_3 \vee X_4 < X_1 \wedge
X_2 ) \\
&&\qquad{}-
P ( Y_1 \wedge Y_3 > Y_2 \vee
Y_4, X_3 \vee X_4 < X_1 \wedge
X_2 )\\
&&\qquad{} -
P ( Y_1 \vee Y_3 < Y_2 \wedge
Y_4, X_3 \vee X_4 < X_1 \wedge
X_2 )\\
&&\quad =
\sum \bigl( P \bigl( Y < y \vert  X_{\nosymbol} =
x_{\nosymbol}^{ ( 1  )} \bigr) - P \bigl( Y < y \vert X_{\nosymbol}
< x_{\nosymbol}^{ ( 1  )} \wedge x_{\nosymbol}^{ ( 2  )}
\bigr) \bigr) \\
&&\quad\qquad{}\times
\bigl( P \bigl( Y < y \vert  X_{\nosymbol} = x_{\nosymbol
}^{ (
2  )}
\bigr) - P \bigl( Y < y \vert X_{\nosymbol} < x_{\nosymbol}^{ ( 1  )}
\wedge x_{\nosymbol
}^{ ( 2
)} \bigr) \bigr) \\
&&\quad\qquad{}\times
P \bigl( Y = y \vert X_{\nosymbol} < x_{\nosymbol}^{ ( 1
)} \wedge
 x_{\nosymbol}^{ ( 2  )} \bigr) \bigl( P \bigl( X <
x^{ ( 1  )} \wedge x^{ ( 2  )} \bigr) \bigr)^2\\
&&\qquad{} +
\sum \bigl( P \bigl( Y < y \vert  X_{\nosymbol} =
x_{\nosymbol}^{ ( 1  )} \bigr) - P \bigl( Y < y \vert X_{\nosymbol}
< x_{\nosymbol}^{ ( 1  )} \wedge x^{ ( 2
)}_{\nosymbol}
\bigr) \bigr)\\
&&\hspace*{5pt}\qquad\qquad{} \times
\bigl( P \bigl( Y < y \vert  X_{\nosymbol} = x_{\nosymbol
}^{ (
2  )}
\bigr) - P \bigl( Y < y \vert X_{\nosymbol} < x_{\nosymbol}^{ ( 1  )}
\wedge x_{\nosymbol
}^{ ( 2
)} \bigr) \bigr) \\
&&\hspace*{5pt}\qquad\qquad{}\times
P \bigl( Y = y \vert X_{\nosymbol} < x_{\nosymbol}^{ ( 1
)} \wedge
 x_{\nosymbol}^{ ( 2  )} \bigr) P \bigl( Y > y \vert
X_{\nosymbol} < x_{\nosymbol}^{ ( 1  )} \wedge
x_{\nosymbol}^{ ( 2  )} \bigr) \bigl( P \bigl( X < x^{ ( 1
)}
\wedge x^{ ( 2  )} \bigr) \bigr)^2\\
&&\qquad{} +
\sum \bigl( P \bigl( Y \leq y \vert  X_{\nosymbol} =
x_{\nosymbol}^{ ( 1  )} \bigr) - P \bigl( Y \leq y \vert
X_{\nosymbol} < x_{\nosymbol}^{ ( 1  )} \wedge
x_{\nosymbol}^{ ( 2  )} \bigr) \bigr)\\
&&\hspace*{5pt}\qquad\qquad{} \times
\bigl( P \bigl( Y \leq y \vert  X_{\nosymbol} =
x_{\nosymbol}^{ ( 2  )} \bigr) - P \bigl( Y \leq y \vert
X_{\nosymbol} < x_{\nosymbol}^{ ( 1  )} \wedge
x^{ ( 2
)}_{\nosymbol} \bigr) \bigr)\\
&&\hspace*{5pt}\qquad\qquad{} \times
P \bigl( Y = y \vert X_{\nosymbol} < x_{\nosymbol}^{ ( 1
)} \wedge
 x_{\nosymbol}^{ ( 2  )} \bigr) \bigl( P \bigl( X <
x^{ ( 1  )} \wedge x^{ ( 2  )} \bigr) \bigr)^2 \\
&&\qquad{}+
\sum \bigl( P \bigl( Y \leq y \vert  X_{\nosymbol} =
x_{\nosymbol}^{ ( 1  )} \bigr) - P \bigl( Y \leq y \vert
X_{\nosymbol} < x_{\nosymbol}^{ ( 1  )} \wedge
x^{ ( 2
)}_{\nosymbol} \bigr) \bigr)\\
&&\hspace*{5pt}\qquad\qquad{} \times
\bigl( P \bigl( Y \leq y \vert  X_{\nosymbol} =
x_{\nosymbol}^{ ( 2  )} \bigr) - P \bigl( Y \leq y \vert
X_{\nosymbol} < x_{\nosymbol}^{ ( 1  )} \wedge
x_{\nosymbol}^{ ( 2  )} \bigr) \bigr)\\
&&\hspace*{5pt}\qquad\qquad {}\times
P \bigl( Y = y \llvert X_{\nosymbol} < x_{\nosymbol}^{ ( 1
)} \wedge
 x_{\nosymbol}^{ ( 2  )} \bigr) P \bigl( Y < y \vert
X_{\nosymbol} < x_{\nosymbol}^{ ( 1  )} \wedge
x_{\nosymbol}^{ ( 2  )} \bigr) \bigl( P \bigl( X < x^{ ( 1
)}
\wedge x^{ ( 2  )} \bigr) \bigr)^2 \\
&&\qquad{}+
\sum \bigl( P \bigl( Y = y \vert  X_{\nosymbol} =
x_{\nosymbol}^{ ( 1  )} \bigr) - P \bigl( Y = y \vert X_{\nosymbol}
< x_{\nosymbol}^{ ( 1  )} \wedge x^{ ( 2  )}
_{\nosymbol} \bigr) \bigr)\\
&&\hspace*{5pt}\qquad\qquad{} \times
\bigl( P \bigl( Y = y \vert  X_{\nosymbol} = x_{\nosymbol
}^{ (
2  )}
\bigr) - P \bigl( Y = y \vert X_{\nosymbol} < x_{\nosymbol}^{ ( 1  )}
\wedge x_{\nosymbol
}^{ ( 2
)} \bigr) \bigr) \\
&&\hspace*{5pt}\qquad\qquad{}\times
P \bigl( Y < y \vert X_{\nosymbol} < x_{\nosymbol}^{ ( 1
)} \wedge
 x_{\nosymbol}^{ ( 2  )} \bigr) P \bigl( Y > y \vert
X_{\nosymbol} < x_{\nosymbol}^{ ( 1  )} \wedge
x_{\nosymbol}^{ ( 2  )} \bigr) \bigl( P \bigl( X < x^{ ( 1
)}
\wedge x^{ ( 2  )} \bigr) \bigr)^2 \\
&&\qquad{}+
3\int \bigl( P \bigl( Y < y \vert  X_{\nosymbol} =
x_{\nosymbol}^{ ( 1  )} \bigr) - P \bigl( Y < y \vert X_{\nosymbol}
< x_{\nosymbol}^{ ( 1  )} \wedge x_{\nosymbol}^{ ( 2  )}
\bigr) \bigr)\\
&&\hspace*{5pt}\qquad\qquad{} \times
\bigl( P \bigl( Y < y \vert  X_{\nosymbol} = x_{\nosymbol
}^{ (
2  )}
\bigr) - P \bigl( Y < y \vert X_{\nosymbol} < x_{\nosymbol}^{ ( 1  )}
\wedge x_{\nosymbol
}^{ ( 2
)} \bigr) \bigr)\\
&&\hspace*{5pt}\qquad\qquad{} \times
P \bigl( Y \in\,\tmop{d}y \vert X_{\nosymbol} < x_{\nosymbol
}^{ ( 1
)}
\wedge x_{\nosymbol}^{ ( 2  )} \bigr) \bigl( P \bigl( X <
x^{ ( 1  )} \wedge x^{ ( 2  )} \bigr) \bigr)^2 .
\end{eqnarray*}
All of the above terms lead to non-negative expressions because of
Lemma~\ref{lem-ins4} (for the first, third and sixth term we take $C =
\Omega$, the set of all possible outcomes). We then see that the expression
can be zero iff $X$ and $Y$ are independent. The condition stated in
Theorem \ref{th1}
is needed to avoid complications when integrating over $x^{ ( 1
)}$
in the application of Lemma \ref{lem-ins4} to terms such as $\sum P  ( Y = y
\vert
 X_{\nosymbol} = x_{\nosymbol}^{ ( 1  )}  )$.
\end{pf*}

\section*{Acknowledgements}
We would like to thank the anonymous referee for useful comments. We would
also like to thank the Associate Editor for many insightful and important
suggestions that greatly improved this paper.

\begin{supplement}
\stitle{A shorter proof of the main theorem for the continuous case and
some miscellaneous further results}
\slink[doi]{10.3150/13-BEJ514SUPP} 
\sdatatype{.pdf}
\sfilename{BEJ514\_supp.pdf}
\sdescription{The supplement contains the following results: (i) a
shorter proof of the main theorem, but only for the continuous case,
(ii) the \cvm\ test as a special case, (iii) a shorter proof of main
theorem for the case that one of the variables is binary, and (iv) a
result for an extension to the case of variables in metric spaces.}
\end{supplement}

%
%

%

\printhistory

\end{document}